\documentclass[12pt,reqno,a4wide]{amsart}

\usepackage{tikz} 
\usepackage{calrsfs}
\usepackage{mathrsfs}

\usepackage{array}
\usepackage{hyperref}

\usetikzlibrary{decorations.pathreplacing}
\usetikzlibrary{fit}

\usepackage{pgfplots}

\allowdisplaybreaks

\catcode`,\active

\catcode`\,12

\begin{document}
	\bibliographystyle{plain}
	
	%
	%
	
	\title
	{Letters of a given type in Catalan words: A continued fraction approach}

	\author[H. Prodinger ]{Helmut Prodinger }
	\address{Department of Mathematics, University of Stellenbosch 7602, Stellenbosch, South Africa
		and
		NITheCS (National Institute for
		Theoretical and Computational Sciences), South Africa.}
	\email{warrenham33@gmail.com}


	\begin{abstract}
Generating functions related to Catalan words and frequencies of digits are obtained using continued fractions. This is fast, elegant, and flexible.
It follows the philosophy of Philippe Flajolet from 1980.
		
	\end{abstract}

	\maketitle

\section{Catalan words}

Dyck paths consist of $n$ up steps and $n$ down steps. They begin at the $x$-axis and end at the $x$-axis and never go below the $x$-axis. Stanley's book \cite{Stanley} contains them and many equivalent combinatorial objects. They might be described by a graph (automaton) checking the conditions.

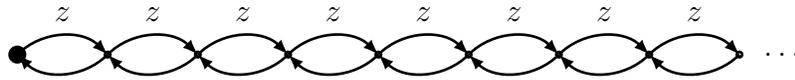
\begin{figure}[h]\label{pso}
	\begin{tikzpicture}[line width=1pt,scale=1.2]
		
		\draw (0,0) circle (0.05cm);
		\fill (0,0) circle (0.1cm);
		\draw (8.5,0.0)node  {$\dots$};

		\foreach \t in{0,1,2,3,4,5,6,7}
		{\draw[-latex](\t,0) to [out=45,in=135](\t+1,0);
			\draw[latex-](\t,0) to [in=-135,out=-45](\t+1,0);
			\draw (\t+1,0)circle (0.03cm);
		}
		\foreach \t in{1,2,3,4,5,6,7,8}
		{\node[label=${z}$]        at (\t-0.5,0.11) {};}

	\end{tikzpicture}
	\caption{Starting and ending at the bullet node. If $2n$ steps are done, the coefficient of $z^n$ is the Catalan number $\frac1{n+1}\binom{2n}{n}$. }
\end{figure}

For the generating function $C(z)$, a decomposition whenever one returns to the beginning leads to $C(z)=\frac1{1-zC(z)}$ and the usual
\begin{equation*}
C(z)=\frac{1-\sqrt{1-4z}}{2z}=\sum_{n\ge0}\frac1{n+1}\binom{2n}{n}z^n.
\end{equation*}
The recursion can be iterated:
\begin{equation*}
C(z)=\cfrac1{1-\cfrac z{1-\cfrac z{1-\cdots}}}
\end{equation*}
For the relation between lattice paths and continued fractions we cite the pioneering paper \cite{Flajolet-CF}.
\begin{figure}[h]\label{pso2}
	\begin{tikzpicture}[line width=1pt,scale=1.2]
		
		\draw (0,0) circle (0.05cm);
		\fill (0,0) circle (0.1cm);
		\draw (8.5,0.0)node  {$\dots$};

		\foreach \t in{0,1,2,3,4,5,6,7}
		{\draw[-latex](\t,0) to [out=45,in=135](\t+1,0);
			\draw[latex-](\t,0) to [in=-135,out=-45](\t+1,0);
		\draw (\t+1,0)circle (0.03cm);
	}
		\foreach \t in{1,2,3,4,5,6,7,8}
		\node[label=${\t}$]        at (\t-0.5,0.11) {};

	\end{tikzpicture}
\caption{Starting and ending at the bullet node. If $2n$ steps are done, the labels of the $n$ forward steps form a \emph{Catalan word}. For $n=3$, these words are possible: $111$,  $112$,  $121$,  $122$,  $123$.  }
\end{figure}
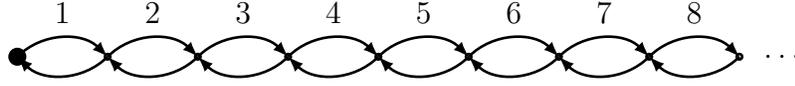

Catalan words are built from letters $1,2,3,\dots$. If the word is written as $a_1a_2\dots a_n$ then $a_1=1$ and $a_{i+1}\le a_i+1$.
I recommend the papers \cite{Callan, Blecher} as a references for Catalan words.

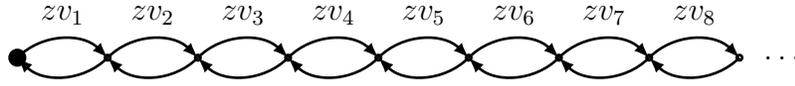
\begin{figure}[h]
	\begin{tikzpicture}[line width=1pt,scale=1.2]
		
		\draw (0,0) circle (0.05cm);
		\fill (0,0) circle (0.1cm);
		\draw (8.5,0.0)node  {$\dots$};

		\foreach \t in{0,1,2,3,4,5,6,7}
		{\draw[-latex](\t,0) to [out=45,in=135](\t+1,0);
			\draw[latex-](\t,0) to [in=-135,out=-45](\t+1,0);
			\draw (\t+1,0)circle (0.03cm);
		}
		\foreach \t in{1,2,3,4,5,6,7,8}
		\node[label=$zv_{\t}$]        at (\t-0.5,0.11) {};

	\end{tikzpicture}
	\caption{Using variables to track forward letters (up steps).}
	\label{pso3}
\end{figure}

If one wants to keep track of what types of forward steps are done, one can use variables $v_1, v_2,\dots$ as in Figure~\ref{pso3}.

\section{Catalan words and continued fractions}

The continued fraction is now
\begin{equation*}
	C(z)=\cfrac1{1-\cfrac{zv_1}{1-\cfrac{zv_2}{1-\ddots}}}
\end{equation*}
\begin{table}[h]
	\caption{$h_i=h_{i-1}-a_ih_{i-2}$,\ $k_i=k_{i-1}-a_ik_{i-2}$}
	\label{la1}
	\begin{tabular}{|l||c|c|c|c|c|c|}
		\hline
		$n$&$-1$& $0 $   & $1$&$2$&$3$&$4$ \\ \hline
		$a_n$ && 1& $zv_1$& $zv_2$    & $zv_3$&$zv_4$ \\
		$h_n$ &  $0$&1&$1$&$1-zv_2$    & $1-zv_2-zv_3$&\dots \\
		$k_n$ &   $1$&1&$1-zv_1$& $1-zv_1-zv_2$    & $1-zv_1-zv_2-zv_3+z^2v_1v_3$&$\dots$ \\
				\hline
	\end{tabular}
	\label{bul}
\end{table}
With the recursive scheme as in Table~\ref{bul}, we find the representation (the traditional method to unwind a continued fraction).
\begin{equation*}
\frac{h_4}{k_4}=\cfrac1{1-\cfrac {zv_1}{1-\cfrac {zv_2}{1-\cfrac {zv_3}{1-zv_4C(z)}}}}
\end{equation*}
and similar in general for $\frac{h_n}{k_n}$. Replacing $C(z)$ by 1 leads to the generating function where forward letters $>n$ do not occur. In other words, these are Dyck paths of bounded height, a very classical subject \cite{deBrKnRi}, and we are in the realm of Chebyshev polynomials. If we replace $C(z)$ by $\frac{1-\sqrt{1-4z}}{2z}$, we have the generating function $F(z;v_1,\dots,v_n)$ and count the semilength of the Catalan words, and the occurrences of letters $1,2,\dots,n$. As an example, we cite
\begin{align*}
\frac{h_5}{k_5}&=\frac{1-zv_5C-v_4z-zv_3+z^2v_3v_5C-zv_2+z^2v_2v_5C+z^2v_2v_4}{\displaystyle{\genfrac{}{}{0pt}{}{1-zv_5C-v_4z-zv_3+z^2v_3v_5C-zv_2+z^2v_2v_5C+z^2v_2v_4-zv_1}{+z^2v_1v_5C+z^2v_1v_4+z^2v_1v_3-z^3v_1v_3v_5C}}}\\
&=1+zv_1+(v_1v_2+v_1^2)z^2+(v_1v_2v_3+v_1v_2^2+2v_1^2v_2+v_1^3)z^3\\&\qquad+(3v_1^2v_2^2+3v_1^3v_2+2v_1^2v_2v_3+2v_1v_2^2v_3+v_1v_2v_3v_4+v_1v_2v_3^2+v_1v_2^3+v_1^4)z^4+\dots;
\end{align*}
\begin{equation*}
\frac{h_3}{k_3}=1+zv_1+(v_1v_2+v_1^2)z^2+\mathbf{(v_1v_2v_3+v_1v_2^2+2v_1^2v_2+v_1^3)}z^3+\dots
\end{equation*}
The 5 terms in bold correspond to $123$, $122$, $112$, $121$, $111$.

Blecher and Knopfmacher \cite{Blecher} were interested in more manageable generating functions, namely by setting all variables $v_j=1$, except for one variable $v_i$, for which we write simply $V$. In an ad hoc notation, we can write these functions for the first few values of $i=1,2,\dots,10$:
\begin{align*}
\mathscr{A}_1(z)&=\frac{1}{1-zVC},\\
\mathscr{A}_2(z)&=\frac{1-zVC}{1-zVC-z},\\
\mathscr{A}_3(z)&=\frac{1-zVC-z}{1-zVC-2z+z^2VC},\\
\mathscr{A}_4(z)&=\frac{1-zVC-2z+z^2VC}{1-zVC-3z+2z^2VC+z^2},\\
\mathscr{A}_5(z)&=\frac{1-zVC-3z+2z^2VC+z^2}{1-zVC-4z+3z^2VC+3z^2-z^3VC},\\
\mathscr{A}_6(z)&=\frac{1-zVC-4z+3z^2VC+3z^2-z^3VC}{1-zVC-5z+4z^2VC+6z^2-3z^3VC-z^3},\\
\mathscr{A}_7(z)&=\frac{1-zVC-5z+4z^2VC+6z^2-3z^3VC-z^3}{1-zVC-6z+5z^2VC+10z^2-6z^3VC-4z^3+z^4VC},\\
\mathscr{A}_8(z)&=\frac{1-zVC-6z+5z^2VC+10z^2-6z^3VC-4z^3+z^4VC}{1-zVC-7z+6z^2VC+15z^2-10z^3VC-10z^3+4z^4VC+z^4},\\
\mathscr{A}_9(z)&=\frac{1-zVC-7z+6z^2VC+15z^2-10z^3VC-10z^3+4z^4VC+z^4}{1-zVC-8z+7z^2VC+21z^2-15z^3VC-20z^3+10z^4VC+5z^4-z^5VC},\\
\mathscr{A}_{10}(z)&=\frac{1-zVC-8z+7z^2VC+21z^2-15z^3VC-20z^3+10z^4VC+5z^4-z^5VC}{1-zVC-9z+8z^2VC+28z^2-21z^3VC-35z^3+20z^4VC+15z^4-5z^5VC-z^5}.
\end{align*}
Expanding, we get for example
\begin{align*}
\mathscr{A}_5(z)&=1+z+2{z}^{2}+5{z}^{3}+14{z}^{4}+ ( 41+V ) {z}^{5}+
( 122+9V+{V}^{2} ) {z}^{6}\\&+ ( 365+52V+11{V}^{2}
+{V}^{3} ) {z}^{7}+ ( 1094+247V+75{V}^{2}+13{V}^{3}+{
	V}^{4} ) {z}^{8}\\&+ ( 3281+1053V +410{V}^{2}+102{V}^{3}+15{V}^
{4}+{V}^{5}) {z}^{9}\\&+ ( 9842+4199V+1975{V}^{2}+629{V
}^{3}+133{V}^{4}+17{V}^{5}+{V}^{6} ) {z}^{10}+\dots
\end{align*}

\section{Counting subwords 12 and variations}

Our goal is to label the occurrences of $12$ with a variable $w$. First, we decompose Catalan words accordingly. We use $d$ for a downstep, the star operator $^*$ as usual in formal languages (arbitrary repetitions) and $\mathcal{C}$ for the set of Catalan words (although with letters $3,4,\dots$ but these letters are not labeled).

An excursion (coming back to the origin for the first time) is then 
\begin{equation*}
1d+12(d2+3\mathcal{C}d)^*dd
\end{equation*}
or in terms of generating functions
\begin{equation*}
	z+z^2w(z+zC(z))^*=	z+\frac{z^2w}{1-z-zC(z)}.
\end{equation*}
A slightly different argument replaces $(d2+3\mathcal{C}d)^*dd$ by $\mathcal{C}d\mathcal{C}d$.
Arbitrary repetitions of this lead to $(z+zw(z+zC(z))^*)^*$, or
\begin{equation*}
F_{12}(z,w)=\frac{1}{ 1-z-{\cfrac {{z}^{2}w}{1-z-{z{ C(z)}}}} }.
\end{equation*}
Note that $F_{12}(z,1)=C(z)$, as nothing is labelled and all Catalan words are considered. A series expansion leads to
\begin{equation*}
	F_{12}(z,w)=1+z+ \left( w+1 \right) {z}^{2}+ \mathbf{\left( 4w+1 \right)} {z}^{3}+
	\left( {w}^{2}+12w+1 \right) {z}^{4}+ \left( 7{w}^{2}+34w+1
	\right) {z}^{5}+\dots
	\end{equation*}
Considering the 5 Catalan words $123$, $122$, $112$, $121$, $111$ of semi-length $3$, we see that $4$ contain the subword $12$ once each time, and $111$ does not contain it.
This explains the bold factor $\mathbf{\left( 4w+1 \right)}$.

The idea of continued fractions allows to count ocurrences of $23$ via 
\begin{align*}
F_{23}(z,w)&=\frac{1}{1-zF_{12}(z,w)}\\
&=1+z+2z^2+\textbf{(w+4)}z^3+(8+6w)z^4+(25w+16+w^2)z^5+ \cdots
\end{align*}
and $\textbf{(w+4)}$ refers to one occurrence of $23$, namely in $123$.

The idea can be iterated as 
\begin{align*}
	F_{34}(z,w)&=\frac{1}{1-zF_{23}(z,w)}=
	1+z+2z^2+5z^3+(13+w)z^4+(34+8w)z^5+\cdots.
\end{align*}
etc. Indeed, $F_{m-1,m}(z,w)$ can be expressed $F_{12}(z,w)$ as a continued fraction.

Now we show how to count occurrences of 123. We start again with excursions that return to the $x$-axis for the first time:
\begin{equation*}
1d+12d\mathcal{C}d+123\mathcal{C}d\mathcal{C}d\mathcal{C}d,
\end{equation*}
in terms of general functions 
\begin{equation*}
z+z^2C+z^3wC^3.
\end{equation*}
Eventually we need arbitrary repetitions of it:
\begin{align*}
F_{123}(z,w)&=\frac1{1-z-z^2C-z^3wC^3}\\
&=1+z+2z^2+(w+4)z^3+(9+5w)z^4+(22+20w)z^5+\cdots.
\end{align*}
Again, we  have 
\begin{equation*}
F_{234}(z,w)=\frac{1}{1-zF_{123}(z,w)}.
\end{equation*}
The presented principles extend to the counting of subwords 1234: We start with
\begin{equation*}
z+z^2C+z^3C^2+wz^4C^4
\end{equation*}
and get
\begin{align*}
	F_{1234}(z,w)&=\frac1{1-z-z^2C-z^3C^2-wz^4C^4}\\
	&=1+z+2z^2+5z^3+(13+w)z^4+(36+6w)z^5+\cdots.
\end{align*}
One can see the general principle (for $h\ge2$):
\begin{align*}
	F_{12\dots h}(z,w)&=\frac1{1-z-z^2C-\dots - z^{h-1}C^{h-2}-wz^hC^h}.
\end{align*}

\section{Final thoughts}

The present approach with continued fractions is fast, elegant, and flexible. It can be recommended for similar problems.


\end{document}